\documentclass[12pt]{article}
\usepackage{amsmath,amssymb,amsthm}

\newtheorem{theorem}{Theorem}
\newtheorem{proposition}{Proposition}
\newtheorem{lemma}[theorem]{Lemma}

\begin{document}

\title{Multivector Functions of a Real Variable\thanks{%
published: \emph{Advances in Applied Clifford Algebras} \textbf{11}(S3),
69-77 (2001).}}
\author{A. M. Moya$^{1}\thanks{%
e-mail: moya@ime.unicamp.br}$, V. V. Fern\'{a}ndez$^{1}\thanks{%
e-mail: vvf@ime.unicamp.br}$ and W. A. Rodrigues Jr.$^{1,2}\thanks{%
e-mail: walrod@ime.unicamp.br or walrod@mpc.com.br}$ \\
$\hspace{-0.5cm}^{1}$ Institute of Mathematics, Statistics and Scientific
Computation\\
IMECC-UNICAMP CP 6065\\
13083-970 Campinas-SP, Brazil\\
$^{2}$ Department of Mathematical Sciences, University of Liverpool\\
Liverpool, L69 3BX, UK}
\date{11/26/2001}
\maketitle

\begin{abstract}
This paper is an introduction to the theory of multivector functions of a
real variable. The notions of limit, continuity and derivative for these
objects are given. The theory of multivector functions of a real variable,
even being similar to the usual theory of vector functions of a real
variable, has some subtle issues which make its presentation worhtwhile.We
refer in particular to the derivative rules involving exterior and Clifford
products, and also to the rule for derivation of a composition of an
ordinary scalar function with a multivector function of a real variable.
\end{abstract}

\tableofcontents

\section{Introduction}

This is paper V of a series of seven. Here, we develop a theory of
multivector functions of a real variable following analogous steps to the
elementary theory of vector functions of a real variable. We introduce the
notions of limit and continuity, and study the concept of derivative. There
are \emph{subtle} points that are emphasized whose understanding is crucial
for the development of a theory of multivector functions of a multivector
variable (as constructed in paper VI of the present series). We give the
complete proofs for the derivative rules involving all the suitable products
of multivector functions of a real variable, and for the composition of an
ordinary scalar function with a multivector function of a real variable.

\section{Multivector Functions of a Real Variable}

Any mapping which takes real numbers of $S\subseteq\mathbb{R}$ to
multivectors of $\bigwedge V$ will be called a \emph{multivector function of
a real variable over }$V.$ In particular, $X:S\rightarrow\bigwedge^{p}V$ is
said to be a $p$\emph{-vector function of a real variable.} And, the special
cases $p=0,$ $p=1,$ $p=2,\ldots,$ etc. are named as a \emph{scalar, vector,
bivector,}$\ldots,$\emph{\ etc. function of a real variable, }respectively.

\subsection{Limit Notion}

We begin by recalling the fundamental concept of $\delta$-neighborhood for a
real number $\lambda_{0}$.

Take any real $\delta>0.$ The set\footnote{%
The symbol $\left| \left. {}\right. \right| $ denotes as usual the \emph{%
absolute value} (or, \emph{module}) function.} $N_{\lambda_{0}}(\delta)=\{%
\lambda\in\mathbb{R}/$ $\left| \lambda-\lambda_{0}\right| <\delta\}, $
clearly a subset of $\mathbb{R},$ is usually called a\emph{\ }$\delta$\emph{%
-neighborhood of }$\lambda_{0}.$ The set $N_{\lambda_{0}}^{\prime}(%
\delta)=N_{\lambda_{0}}(\delta)-\{\lambda_{0}\},$ i.e., $N_{\lambda_{0}}^{%
\prime}(\delta )=\{\lambda\in\mathbb{R}/$ $0<\left|
\lambda-\lambda_{0}\right| <\delta\},$ is said to be a \emph{reduced }$\delta
$\emph{-neighborhood of }$\lambda_{0}.$

We recall now the important concept of cluster point and interior point of $%
S\subseteq\mathbb{R}.$

A real number $\lambda_{0}$ is said to be a \emph{cluster point of }$S$ if
and only if for every $N_{\lambda_{0}}(\delta):N_{\lambda_{0}}^{\prime}(%
\delta)\cap S\neq\emptyset,$ i.e., all reduced $\delta$-neighborhood of $%
\lambda_{0}$ contains at least one real number of $S.$

A real number $\lambda_{0}$ is said to be a \emph{interior point of }$S$%
\emph{\ }if and only if there exists $N_{\lambda_{0}}(\delta)$ such that $%
N_{\lambda_{0}}(\delta)\subseteq S,$ i.e., any real number of some $\delta $%
-neighborhood of $\lambda_{0}$ belongs also to $S$.

Note that all interior point of $S$ is also cluster point of $S.$

If the set of interior point of $S$ coincides with $S,$ i.e., all real
number of $S$ is also interior point of $S,$ then $S$ is said to be an \emph{%
open subset of }$\mathbb{R}.$

Next we introduce the concept of \emph{norm} of a multivector $X.$

Assume that $\bigwedge V$ has been endowed with an \emph{euclidean scalar
product }$(\cdot)$, as e.g., by taking any fixed basis $\{b_{k}\}$ for $V$
and its dual basis $\{\beta^{k}\}$ for $V^{*},$ etc. See paper I of this
series \cite{1}. As we already know, the euclidean scalar product is always
\emph{definite positive, }i.e., for all $X\in\bigwedge V:X\cdot X\geq0$ and $%
X\cdot X=0$ if and only if $X=0.$

This property of the scalar product permit us to introduce the \emph{norm of
a multivector }$X$ as being the non-negative real number $\left\| X\right\| $
given by
\begin{equation}
\left\| X\right\| =\sqrt{X\cdot X}.   \label{5.0a}
\end{equation}
We read $\left\| X\right\| $ as the \emph{norm of} $X.$

The norm of multivectors satisfies the following two usual inequalities:

\textbf{n1} The Cauchy-Schwarz inequality, i.e., for all $X,Y\in\bigwedge V$
\begin{equation}
\left| X\cdot Y\right| \leq\left\| X\right\| \left\| Y\right\| .
\label{5.0b}
\end{equation}

\textbf{n2} The triangular inequality, i.e., for all $X,Y\in\bigwedge V$
\begin{equation}
\left\| X+Y\right\| \leq\left\| X\right\| +\left\| Y\right\| .   \label{5.0c}
\end{equation}

The first inequality follows from the fact that $(\cdot)$ is positive
definite. The second one is an immediate consequence of the first one.

Take $S\subseteq\mathbb{R}.$ Let $X:S\rightarrow\bigwedge V$ be any
multivector function of a real variable and take $\lambda_{0}\in S$ to be a
cluster point of $S.$

A multivector $L$ is said to be the \emph{limit of }$X(\lambda)$\emph{\ for }%
$\lambda$\emph{\ approaching to }$\lambda_{0}$ if and only if for every real
$\varepsilon>0$ there exists some real $\delta>0$ such that if for all $%
\lambda\in S$ and $0<\left| \lambda-\lambda_{0}\right| <\delta,$ then $%
\left\| X(\lambda)-L\right\| <\varepsilon.$ It is denoted by $\underset{%
\lambda\rightarrow\lambda_{0}}{\lim}X(\lambda)=L.$

In particular, a scalar function of a real variable is just an ordinary real
function and, as we can see, the above definition of limit is reduced to the
ordinary definition of limit which appears in real analysis.

\begin{proposition}
Let $X:S\rightarrow\bigwedge V$ and $Y:S\rightarrow\bigwedge V$ be two
multivector functions of a real variable. If there exist $\underset{%
\lambda\rightarrow\lambda_{0}}{\lim}X(\lambda)$ and $\underset{\lambda
\rightarrow\lambda_{0}}{\lim}Y(\lambda),$ then there exists $\underset{%
\lambda\rightarrow\lambda_{0}}{\lim}(X+Y)(\lambda)$ and
\begin{equation}
\underset{\lambda\rightarrow\lambda_{0}}{\lim}(X+Y)(\lambda)=\underset{%
\lambda\rightarrow\lambda_{0}}{\lim}X(\lambda)+\underset{\lambda
\rightarrow\lambda_{0}}{\lim}Y(\lambda).   \label{5.1a}
\end{equation}
\end{proposition}

\begin{proof}
Let $\underset{\lambda\rightarrow\lambda_{0}}{\lim}X(\lambda)=L_{1}$ and
$\underset{\lambda\rightarrow\lambda_{0}}{\lim}Y(\lambda)=L_{2}.$ Then, we
must prove that $\underset{\lambda\rightarrow\lambda_{0}}{\lim}(X+Y)(\lambda
)=L_{1}+L_{2}.$

Given an arbitrary real $\varepsilon>0,$ since $\underset{\lambda
\rightarrow\lambda_{0}}{\lim}X(\lambda)=L_{1}$ and $\underset{\lambda
\rightarrow\lambda_{0}}{\lim}Y(\lambda)=L_{2},$ there are two real numbers
$\delta_{1}>0$ and $\delta_{2}>0$ such that
\begin{align*}
\left\|  X(\lambda)-L_{1}\right\|   & <\frac{\varepsilon}{2},\text{ for
}\lambda\in S\text{ and }0<\left|  \lambda-\lambda_{0}\right|  <\delta_{1},\\
\left\|  Y(\lambda)-L_{2}\right\|   & <\frac{\varepsilon}{2},\text{ for
}\lambda\in S\text{ and }0<\left|  \lambda-\lambda_{0}\right|  <\delta_{2}.
\end{align*}

Thus, there is a real $\delta=\min\{\delta_{1},\delta_{2}\}$ such that
\[
\left\|  X(\lambda)-L_{1}\right\|  <\frac{\varepsilon}{2}\text{ and }\left\|
Y(\lambda)-L_{2}\right\|  <\frac{\varepsilon}{2},\text{ }
\]
for $\lambda\in S$ and $0<\left|  \lambda-\lambda_{0}\right|  <\delta.$ Hence,
by using eq.(\ref{5.0c}) it follows that
\begin{align*}
\left\|  (X+Y)(\lambda)-(L_{1}+L_{2})\right\|   & =\left\|  X(\lambda
)-L_{1}+Y(\lambda)-L_{2}\right\| \\
& \leq\left\|  X(\lambda)-L_{1}\right\|  +\left\|  Y(\lambda)-L_{2}\right\| \\
& <\frac{\varepsilon}{2}+\frac{\varepsilon}{2}=\varepsilon,
\end{align*}
for $\lambda\in S$ and $0<\left|  \lambda-\lambda_{0}\right|  <\delta.$

Therefore, for any $\varepsilon>0$ there exists a $\delta>0$ such that if
$\lambda\in S$ and $0<\left|  \lambda-\lambda_{0}\right|  <\delta,$ then
$\left\|  (X+Y)(\lambda)-(L_{1}+L_{2})\right\|  <\varepsilon.$
\end{proof}

\begin{proposition}
Let $\phi:S\rightarrow\mathbb{R}$ and $X:S\rightarrow\bigwedge V$ be an
ordinary real function and a multivector function of a real variable. If
there exist $\underset{\lambda\rightarrow\lambda_{0}}{\lim}\phi(\lambda)$
(the ordinary limit) and $\underset{\lambda\rightarrow\lambda_{0}}{\lim}%
X(\lambda),$ then there exists $\underset{\lambda\rightarrow\lambda_{0}}{%
\lim }(\phi X)(\lambda)$ and
\begin{equation}
\underset{\lambda\rightarrow\lambda_{0}}{\lim}(\phi X)(\lambda)=\underset{%
\lambda\rightarrow\lambda_{0}}{\lim}\phi(\lambda)\underset{\lambda
\rightarrow\lambda_{0}}{\lim}X(\lambda).   \label{5.1b}
\end{equation}
\end{proposition}

\begin{proof}
Let $\underset{\lambda\rightarrow\lambda_{0}}{\lim}\phi(\lambda)=\phi_{0}$ and
$\underset{\lambda\rightarrow\lambda_{0}}{\lim}X(\lambda)=X_{0}.$ Then, we
must prove that $\underset{\lambda\rightarrow\lambda_{0}}{\lim}(\phi
X)(\lambda)=\phi_{0}X_{0}.$

First, since $\underset{\lambda\rightarrow\lambda_{0}}{\lim}\phi(\lambda
)=\phi_{0}$ it can be found a $\delta_{1}>0$ such that
\[
\left|  \phi(\lambda)-\phi_{0}\right|  <1,\text{ whenever }\lambda\in S\text{
and }0<\left|  \lambda-\lambda_{0}\right|  <\delta_{1},
\]
i.e.,
\[
\left|  \phi(\lambda)\right|  <1+\left|  \phi_{0}\right|  ,\text{ whenever
}\lambda\in S\text{ and }0<\left|  \lambda-\lambda_{0}\right|  <\delta_{1}.
\]
Where the triangular inequality for real numbers $\left|  \alpha\right|
-\left|  \beta\right|  \leq\left|  \alpha-\beta\right|  $ was used.

Now, taken an arbitrary $\varepsilon>0,$ since $\underset{\lambda
\rightarrow\lambda_{0}}{\lim}\phi(\lambda)=\phi_{0}$ and $\underset
{\lambda\rightarrow\lambda_{0}}{\lim}X(\lambda)=X_{0},$ they can be found a
$\delta_{2}>0$ and a $\delta_{3}>0$ such that
\begin{align*}
\left|  \phi(\lambda)-\phi_{0}\right|   & <\frac{\varepsilon}{2(1+\left\|
X_{0}\right\|  )},\text{ whenever }\lambda\in S\text{ and }0<\left|
\lambda-\lambda_{0}\right|  <\delta_{2},\\
\left\|  X(\lambda)-X(\lambda_{0})\right\|   & <\frac{\varepsilon}{2(1+\left|
\phi_{0}\right|  )},\text{ whenever }\lambda\in S\text{ and }0<\left|
\lambda-\lambda_{0}\right|  <\delta_{3}.
\end{align*}

Thus, given an arbitrary $\varepsilon>0$ there is a $\delta=\min\{\delta
_{1},\delta_{2},\delta_{3}\}$ such that
\begin{align*}
\left|  \phi(\lambda)\right|   & <1+\left|  \phi_{0}\right|  ,\\
\left|  \phi(\lambda)-\phi_{0}\right|   & <\frac{\varepsilon}{2(1+\left\|
X_{0}\right\|  )},\\
\left\|  X(\lambda)-X_{0})\right\|   & <\frac{\varepsilon}{2(1+\left|
\phi_{0}\right|  )},
\end{align*}
for $\lambda\in S$ and $0<\left|  \lambda-\lambda_{0}\right|  <\delta.$ Hence,
it follows that
\begin{align*}
\left\|  (\phi X)(\lambda)-\phi_{0}X_{0}\right\|   & =\left\|  \phi
(\lambda)(X(\lambda)-X_{0})+(\phi(\lambda)-\phi_{0})X_{0}\right\| \\
& \leq\left|  \phi(\lambda)\right|  \left\|  X(\lambda)-X_{0}\right\|
+\left|  \phi(\lambda)-\phi_{0}\right|  \left\|  X_{0}\right\| \\
& <\left|  \phi(\lambda)\right|  \left\|  X(\lambda)-X_{0}\right\|  +\left|
\phi(\lambda)-\phi_{0}\right|  (1+\left\|  X_{0}\right\|  )\\
& <(1+\left|  \phi_{0}\right|  )\frac{\varepsilon}{2(1+\left|  \phi
_{0}\right|  )}+\frac{\varepsilon}{2(1+\left\|  X_{0}\right\|  )}(1+\left\|
X_{0}\right\|  )=\varepsilon,
\end{align*}
for $\lambda\in S$ and $0<\left|  \lambda-\lambda_{0}\right|  <\delta$. In the
proof above we use some properties of the norm of multivectors.

Therefore, for any $\varepsilon>0$ there exists a $\delta>0$ such that if
$\lambda\in S$ and $0<\left|  \lambda-\lambda_{0}\right|  <\delta,$ then
$\left\|  (\phi X)(\lambda)-\phi_{0}X_{0}\right\|  <\varepsilon.$
\end{proof}

\begin{lemma}
There exists $\underset{\lambda\rightarrow\lambda_{0}}{\lim}X(\lambda) $ if
and only if there exist any one of the ordinary limits, either $\underset{%
\lambda\rightarrow\lambda_{0}}{\lim}X^{J}(\lambda)$ or $\underset{%
\lambda\rightarrow\lambda_{0}}{\lim}X_{J}(\lambda).$ It holds
\begin{equation}
\underset{\lambda\rightarrow\lambda_{0}}{\lim}X(\lambda)=\underset{J}{\sum }%
\frac{1}{\nu(J)!}\underset{\lambda\rightarrow\lambda_{0}}{\lim}%
X^{J}(\lambda)e_{J}=\underset{J}{\sum}\frac{1}{\nu(J)!}\underset{\lambda
\rightarrow\lambda_{0}}{\lim}X_{J}(\lambda)e^{J}.   \label{5.1c}
\end{equation}
\end{lemma}

\begin{proof}
It is an immediate consequence of eqs.(\ref{5.1a}) and (\ref{5.1b}).
\end{proof}

\begin{proposition}
Let $X:S\rightarrow\bigwedge V$ and $Y:S\rightarrow\bigwedge V$ be two
multivector functions of a real variable. We can define the products $X\ast
Y:S\rightarrow\bigwedge V$ such that $(X\ast Y)(\lambda)=X(\lambda)\ast
Y(\lambda)$ where $\ast$ holds for either $(\wedge),$ $(\cdot),$ $%
(\lrcorner\llcorner)$ or $($\emph{Clifford product}$).$ If there exist $%
\underset{\lambda\rightarrow\lambda_{0}}{\lim}X(\lambda)$ and $\underset{%
\lambda\rightarrow\lambda_{0}}{\lim}Y(\lambda),$ then there exists $\underset%
{\lambda\rightarrow\lambda_{0}}{\lim}(X\ast Y)(\lambda)$ and
\begin{equation}
\underset{\lambda\rightarrow\lambda_{0}}{\lim}(X\ast Y)(\lambda)=\underset{%
\lambda\rightarrow\lambda_{0}}{\lim}X(\lambda)\ast\underset{\lambda
\rightarrow\lambda_{0}}{\lim}Y(\lambda).   \label{5.1d}
\end{equation}
\end{proposition}

\begin{proof}
It is an immediate consequence of eq.(\ref{5.1c}).
\end{proof}

\subsection{Continuity Notion}

Take $S\subseteq\mathbb{R}.$ A multivector function of a real variable $%
X:S\rightarrow\bigwedge V$ is said to be \emph{continuous at }$%
\lambda_{0}\in S$ if and only if there exists\footnote{%
See that $\lambda_{0}$ has to be cluster point of $S.$} $\underset{%
\lambda\rightarrow\lambda_{0}}{\lim }X(\lambda)$ and
\begin{equation}
\underset{\lambda\rightarrow\lambda_{0}}{\lim}X(\lambda)=X(\lambda _{0}).
\label{5.2a}
\end{equation}

\begin{lemma}
The multivector function $\lambda\mapsto X(\lambda)$ is continuous at $%
\lambda_{0}$ if and only if whichever component scalar function either $%
\lambda\mapsto X^{J}(\lambda)$ or $\lambda\mapsto X_{J}(\lambda)$ is
continuous at $\lambda_{0}.$
\end{lemma}

\begin{proposition}
Let $X:S\rightarrow\bigwedge V$ and $Y:S\rightarrow\bigwedge V$ be two
continuous functions at $\lambda_{0}\in S.$

The addition $X+Y:S\rightarrow\bigwedge V$ such that $(X+Y)(\lambda
)=X(\lambda)+Y(\lambda)$ and the products $X\ast Y:S\rightarrow\bigwedge V $
such that $(X\ast Y)(\lambda)=X(\lambda)\ast Y(\lambda),$ where $\ast$ means
either $(\wedge),$ $(\cdot),$ $(\lrcorner\llcorner)$ or $($\emph{Clifford
product}$),$ are also continuous functions at $\lambda_{0}.$
\end{proposition}

\begin{proof}
 It is an immediate consequence of eqs.(\ref{5.1a}) and
(\ref{5.1d}).
\end{proof}

\begin{proposition}
Let $\phi:S\rightarrow\mathbb{R}$ and $X:\mathbb{R}\rightarrow\bigwedge V$
be two continuous functions, the first one at $\lambda_{0}\in S$ and the
second one at $\phi(\lambda_{0})\in R.$

The composition $X\circ\phi:S\rightarrow\bigwedge V$ such that $X\circ
\phi(\lambda)=X(\phi(\lambda))$ is a continuous function at $\lambda_{0}.$
\end{proposition}

\subsection{Derivative}

Take $S\subseteq\mathbb{R}$ be an open set of $\mathbb{R}$. A multivector
function of a real variable $X:S\rightarrow\bigwedge V$ is said to be \emph{%
derivable at }$\lambda_{0}\in S$ if and only if there exists $\underset{%
\lambda\rightarrow\lambda_{0}}{\lim}\dfrac{X(\lambda)-X(\lambda _{0})}{%
\lambda-\lambda_{0}}.$ This \emph{multivector-limit} is usually called the
\emph{derivative of }$X$\emph{\ at }$\lambda_{0}\in S,$\emph{\ }and often
denoted by $X^{\prime}(\lambda_{0}),$ i.e.,
\begin{equation}
X^{\prime}(\lambda_{0})=\underset{\lambda\rightarrow\lambda_{0}}{\lim}\dfrac{%
X(\lambda)-X(\lambda_{0})}{\lambda-\lambda_{0}}.   \label{5.3a}
\end{equation}

So that, the derivability of $X$ at $\lambda_{0}$ means the existence of
derivative of $X$ at $\lambda_{0}.$

\begin{lemma}
Associated to any multivector function $X,$ derivable at $\lambda_{0},$
there exists a multivector function $\xi_{\lambda_{0}},$ continuous at $%
\lambda _{0},$ such that
\begin{equation}
\xi_{\lambda_{0}}(\lambda_{0})=0   \label{5.3b}
\end{equation}
and for all $\lambda\in S$ it holds
\begin{equation}
X(\lambda)=X(\lambda_{0})+(\lambda-\lambda_{0})X^{\prime}(\lambda
_{0})+(\lambda-\lambda_{0})\xi_{\lambda_{0}}(\lambda).   \label{5.3c}
\end{equation}
\end{lemma}

\begin{proof}
Since $X$ is derivable at $\lambda_{0}$ we can define $\xi_{\lambda_{0}}$ by
\[
\xi_{\lambda_{0}}(\lambda)=\left\{
\begin{array}
[c]{cc}%
0 & \text{for }\lambda=\lambda_{0}\\
\dfrac{X(\lambda)-X(\lambda_{0})}{\lambda-\lambda_{0}}-X^{\prime}(\lambda
_{0}) & \text{for }\lambda\neq\lambda_{0}%
\end{array}
\right.  .
\]

We see that $\xi_{\lambda_{0}}(\lambda_{0})=0$ and by taking limit of
$\xi_{\lambda_{0}}(\lambda)$ for $\lambda\rightarrow\lambda_{0}$ we have
\[
\underset{\lambda\rightarrow\lambda_{0}}{\lim}\xi_{\lambda_{0}}(\lambda
)=\underset{\lambda\rightarrow\lambda_{0}}{\lim}(\dfrac{X(\lambda
)-X(\lambda_{0})}{\lambda-\lambda_{0}}-X^{\prime}(\lambda_{0}))=X^{\prime
}(\lambda_{0})-X^{\prime}(\lambda_{0})=0.
\]
It follows that $\xi_{\lambda_{0}}$ is continuous at $\lambda_{0}$ and so the
first statement holds.

On another way, for $\lambda\neq\lambda_{0}$ we get the multivector identity
\[
X(\lambda)=X(\lambda_{0})+(\lambda-\lambda_{0})X^{\prime}(\lambda
_{0})+(\lambda-\lambda_{0})\xi_{\lambda_{0}}(\lambda)
\]
but, for $\lambda=\lambda_{0}$ it is trivially true. Thus, the second
statement holds.
\end{proof}

As happens in real analysis, derivability implies continuity. Indeed, by
taking limits for $\lambda\rightarrow\lambda_{0}$ on both sides of eq.(\ref%
{5.3c}) we get $\underset{\lambda\rightarrow\lambda_{0}}{\lim }%
X(\lambda)=X(\lambda_{0}).$

\subsubsection{Derivation Rules}

Take two open subset of $\mathbb{R},$ say $S_{1}$ and $S_{2},$ such that $%
S_{1}\cap S_{2}\neq\emptyset.$

\begin{theorem}
Let $S_{1}\ni\lambda\mapsto X(\lambda)\in\bigwedge V$ and $S_{2}\ni
\lambda\mapsto Y(\lambda)\in\bigwedge V$ be two derivable functions at $%
\lambda_{0}\in S_{1}\cap S_{2}.$

The addition $S_{1}\cap S_{2}\ni\lambda\mapsto(X+Y)(\lambda)\in\Lambda V $
such that $(X+Y)(\lambda)=X(\lambda)+Y(\lambda)$ and the products $S_{1}\cap
S_{2}\ni\lambda\mapsto(X\ast Y)(\lambda)\in\Lambda V$ such that $(X\ast
Y)(\lambda)=X(\lambda)\ast Y(\lambda),$ where $\ast$ means either $(\wedge)$%
, $(\cdot),$ $(\lrcorner\llcorner)$ or $($Clifford product$),$ are also
derivable functions at $\lambda_{0}.$

The derivatives of $X+Y$ and $X\ast Y$ at $\lambda_{0}$ are given by
\begin{equation}
(X+Y)^{\prime}(\lambda_{0})=X^{\prime}(\lambda_{0})+Y^{\prime}(\lambda _{0})
\label{5.3d}
\end{equation}
and
\begin{equation}
(X\ast Y)^{\prime}(\lambda_{0})=X^{\prime}(\lambda_{0})\ast Y(\lambda
_{0})+X(\lambda_{0})\ast Y^{\prime}(\lambda_{0}).   \label{5.3e}
\end{equation}
\end{theorem}

\begin{proof}
We only need to verify that
\[
\underset{\lambda\rightarrow\lambda_{0}}{\lim}\frac{(X+Y)(\lambda
)-(X+Y)(\lambda_{0})}{\lambda-\lambda_{0}}=X^{\prime}(\lambda_{0})+Y^{\prime
}(\lambda_{0})
\]
and that
\[
\underset{\lambda\rightarrow\lambda_{0}}{\lim}\frac{(X*Y)(\lambda
)-(X*Y)(\lambda_{0})}{\lambda-\lambda_{0}}=X^{\prime}(\lambda_{0}%
)*Y(\lambda_{0})+X(\lambda_{0})*Y^{\prime}(\lambda_{0}).
\]

First, we set the following multivector identities which hold for all
$\lambda\neq\lambda_{0}$
\[
\frac{(X+Y)(\lambda)-(X+Y)(\lambda_{0})}{\lambda-\lambda_{0}}=\frac
{X(\lambda)-X(\lambda_{0})}{\lambda-\lambda_{0}}+\frac{Y(\lambda
)-Y(\lambda_{0})}{\lambda-\lambda_{0}}
\]
and
\[
\frac{(X*Y)(\lambda)-(X*Y)(\lambda_{0})}{\lambda-\lambda_{0}}=\frac
{X(\lambda)-X(\lambda_{0})}{\lambda-\lambda_{0}}*Y(\lambda_{0})+X(\lambda
)*\frac{Y(\lambda)-Y(\lambda_{0})}{\lambda-\lambda_{0}}.
\]

Now, by taking limits for $\lambda\rightarrow\lambda_{0}$ on both sides of
these multivector identities, using the equation\footnote{We have used the
fact that for $X$, derivability implies in continuity.}: $\underset
{\lambda\rightarrow\lambda_{0}}{\lim}X(\lambda)=X(\lambda_{0}),$ we get the
expected results.
\end{proof}

\begin{theorem}
Let $\phi:S\rightarrow\mathbb{R}$ and $X:\mathbb{R}\rightarrow\bigwedge V$
be two derivable functions, the first one at $\lambda_{0}\in S$ and the
second one at $\phi(\lambda_{0})\in\mathbb{R}.$

The composition $X\circ\phi:S\rightarrow\bigwedge V$ such that $X\circ
\phi(\lambda)=X(\phi(\lambda))$ is a derivable function at $\lambda_{0}$ and
its derivative at $\lambda_{0}$ is given by
\begin{equation}
(X\circ\phi)^{\prime}(\lambda_{0})=\phi^{\prime}(\lambda_{0})X^{\prime}(%
\phi(\lambda_{0})).   \label{5.3f}
\end{equation}
\end{theorem}

\begin{proof}
We must prove that
\[
\underset{\lambda\rightarrow\lambda_{0}}{\lim}\frac{X\circ\phi(\lambda
)-X\circ\phi(\lambda_{0})}{\lambda-\lambda_{0}}=\phi^{\prime}(\lambda
_{0})X^{\prime}(\phi(\lambda_{0})).
\]

Since $X$ is derivable at $\phi(\lambda_{0}),$ there is a multivector function
$\mu\mapsto\xi_{\phi(\lambda_{0})}(\mu),$ continuous at $\phi(\lambda_{0}),$
such that for all $\mu\in\mathbb{R}$
\[
X(\mu)=X(\phi(\lambda_{0}))+(\mu-\phi(\lambda_{0}))X^{\prime}(\phi(\lambda
_{0}))+(\mu-\phi(\lambda_{0}))\xi_{\phi(\lambda_{0})}(\mu),
\]
where $\xi_{\phi(\lambda_{0})}(\phi(\lambda_{0}))=0.$

Now, the following multivector identity (as can be easily shown) holds for all
$\lambda\neq\lambda_{0}$,
\[
\frac{X\circ\phi(\lambda)-X\circ\phi(\lambda_{0})}{\lambda-\lambda_{0}}%
=\frac{\phi(\lambda)-\phi(\lambda_{0})}{\lambda-\lambda_{0}}X^{\prime}%
(\phi(\lambda_{0}))+\frac{\phi(\lambda)-\phi(\lambda_{0})}{\lambda-\lambda
_{0}}\xi_{\phi(\lambda_{0})}\circ\phi(\lambda).
\]

Now, by taking limits for $\lambda\rightarrow\lambda_{0}$ on both sides, using
the equation\footnote{It was used that composition of $\phi$ with $\xi
_{\phi(\lambda_{0})},$ where $\phi$ is continuous at $\lambda_{0}$ and
$\xi_{\phi(\lambda_{0})} $ is continuous at $\phi(\lambda_{0}),$ is continuous
at $\lambda_{0}$.}: $\underset{\lambda\rightarrow\lambda_{0}}{\lim}\xi
_{\phi(\lambda_{0})}\circ\phi(\lambda)=0,$ we get the required result.
\end{proof}

\section{Conclusions}

In this paper we introduced the concept of multivector functions of a real
variable, and the notions of limit and continuity for them, and studied the
concept of derivative of these objects. Although our theory of multivector
functions of a real variable parallels the theory of vector functions of a
real variable, we believe that our presentation is worthwhile, since it
treats some subtle points as, e.g., the derivative rules involving all the
suitable products of the multivector functions of a real variable. The
generalization of these ideas towards a general theory of multivector
functions of several real variables can be done without great difficulty.

The results developed in this paper are essential ingredients for papers VI
and VII of the present series of papers, where we obtain important results
concerning to the theory of multivector functions of a multivector variable,
and to the theory of multivector functionals.

Before ending, we quote that the concept of multivector functions (of real
variable or multivector variable) has been first introduced in \cite{2}, and
used together with the notion of multivector functionals by some authors, in
order to study problems ranging from linear algebra to applications to
physical sciences and engineering (e.g., \cite{3}\cite{4}). We believe that
our approach is a real contribution to those presentations of these subjects.%
\vspace{0.1in}

\textbf{Acknowledgement}: V. V. Fern\'{a}ndez is grateful to FAPESP for a
posdoctoral fellowship. W. A. Rodrigues Jr. is grateful to CNPq for a senior
research fellowship (contract 201560/82-8) and to the Department of
Mathematics of the University of Liverpool for the hospitality. Authors are
also grateful to Drs. P. Lounesto, I. Porteous, and J. Vaz, Jr. for their
interest on our research and useful discussions.

\end{document}